\documentclass{amsart}

\newtheorem{theorem}{Theorem}[section]
\newtheorem{lemma}[theorem]{Lemma}
\newtheorem{corollary}[theorem]{Corollary}

\usepackage[main=english]{babel} 
\usepackage{amsmath,amssymb,amscd,amsfonts,amsthm,amsrefs} 
\usepackage{microtype} 

\usepackage{tikz}
\usetikzlibrary{arrows.meta,decorations.markings}
\tikzset{
	path arrow/.style={
			postaction={decorate},
			decoration={markings,mark=at position .53 with
						{\arrow{Stealth[length=4.3pt,width=3.7pt]}}}},
	reverse path arrow/.style={
			postaction={decorate},
			decoration={markings,mark=at position .53 with
						{\arrowreversed{Stealth[length=4.3pt,width=3.7pt]}}}},
	map arrow/.style={-{Stealth[length=5pt,width=4pt]}},
	vertex/.style={circle,fill=black,inner sep=0pt,minimum size=3.2pt}
}

\DeclareMathOperator{\Cay}{Cay}
\newcommand{\abs}[1]{\ensuremath{\lvert#1\rvert}}
\newcommand{\iso}{\cong}
\newcommand{\FF}{\mathbb{F}}

\begin{document}

\begin{abstract}
	Erdős and Gyárfás conjectured that every graph with minimum degree at least three contains a power-of-two-length cycle. We prove every Cayley graph of odd order satisfies the conjecture. As a corollary, we find that Cayley graphs of finite nilpotent groups satisfy the conjecture as well.
\end{abstract}

\author{Michael C. Burkhart}
\curraddr{University of Chicago, Chicago, Illinois}
\email{burkh4rt@uchicago.edu}

\subjclass[2020]{Primary 05C38, 05C25}


\title{Odd-ordered Cayley graphs satisfy Erd\H{o}s--Gy\'arf\'as}

\maketitle

\section{Introduction}

In 1995, Erdős and Gyárfás~\cite{Erd97} conjectured that if a finite graph has minimum degree at least three, it contains a simple cycle of length $2^n, n \ge 2$. Some partial results are known.
Daniel and Shauger~\cite{Dan01} proved the conjecture for claw-free planar graphs, i.e. those without induced $K_{1,3}$.  Heckman and Krakovski~\cite{Hec13} then proved the conjecture for 3-connected cubic planar graphs. Gao and Shan~\cite{Gao22} showed the conjecture must hold for $P_8$-free graphs in work that has been subsequently extended~\cite{Hu24,Heg24} to $P_{13}$-free graphs. It has been shown~\cite{Mar04,Now14,Tra26,Duc26} that a potential counterexample must contain at least 41 vertices, 50 vertices if cubic, and 67 if bipartite.
In the case of large degree, Liu and Montgomery~\cite{Liu23} proved that there exists $d>0$ such that all graphs with average degree at least $d$ contain a power-of-two-length cycle.

For a group $G$ with symmetric generating set $S\subseteq G \setminus\{1\}$, let $\Gamma=\Cay(G,S)$ denote the simple undirected graph with vertices $G$ and edges $\{g, gs\}_{g\in G,s\in S}$. If $S$ has size two, as can be the case with cyclic groups or the two transpositions that generate $S_3$, $\Gamma$ has minimum degree two and so satisfies the conjecture trivially. Otherwise, for abelian $G$ and arbitrary $S, \abs{S}\ge 3$, $\Gamma$ contains the 4-cycle corresponding to the word $sts^{-1}t^{-1}$ for any $s,t \in S, s\ne t^{\pm 1}$. In 2018, Ghaffari and Mostaghim~\cite{Gha18} showed that the conjecture also holds for Cayley graphs of generalized quaternion, dihedral, and semidihedral groups, as well as groups of order $p^3$ for primes $p$. Ghasemi and Varmazyar~\cite{Gha21} extended these results to Cayley graphs of groups having order $2p^2$ or $4p$ for primes $p$. In this work, we show the conjecture holds for all Cayley graphs of odd order:
\begin{theorem}
	\label{thm:main}
	If $G$ is a group of odd order with symmetric generating set $S, \abs{S}\ge 3$, then $\Cay(G,S)$ contains a power-of-two-length cycle.
\end{theorem}

In a Cayley graph $\Cay(G,S)$, we define a \emph{symmetric figure-eight} to be a pair of power-of-two-length simple cycles that intersect at precisely one vertex and where one cycle is a $G$-translate of the other, up to change of basepoint. As translates of each other, the two cycles will have the same length. Our proof shows that every Cayley graph corresponding to a noncyclic group of odd order contains a symmetric figure-eight. After handling base cases, our argument relies on the ability to lift a symmetric figure-eight from a quotient by a normal abelian subgroup to its cover, as follows:
\begin{lemma}
	\label{lem:main}
	Let $G=\langle S \rangle$ be a group of odd order. For abelian $A\triangleleft G$, let $\overline S = \{sA\}_{s\in S} \setminus \{A\}$. If $\Cay(G/A, \overline S)$ contains a symmetric figure-eight, then so does $\Cay(G,S)$.
\end{lemma}

We can also lift a power-of-two-length cycle from a quotient by a normal 2-subgroup to its cover:
\begin{lemma}
	\label{lem:index2}
	Let $G=\langle S \rangle$ be a group with normal 2-subgroup $P\triangleleft G$. If $\Cay(G/P, \overline S)$ contains a power-of-two-length cycle where $\overline S = \{sP\}_{s\in S}\setminus \{P\}$, then $\Cay(G,S)$ contains a power-of-two-length cycle.
\end{lemma}

This, together with Theorem~\ref{thm:main}, allows us to establish an analogous result for nilpotent groups:
\begin{corollary}
	\label{cor:nilp}
	If $G$ is finite nilpotent with symmetric generating set $S, \abs{S}\ge 3$, then $\Cay(G,S)$ contains a power-of-two-length cycle.
\end{corollary}

In the remainder of this section, we outline some notation. Section~\ref{s:odd} concerns groups of odd order and contains the proofs of Lemma~\ref{lem:main} and Theorem~\ref{thm:main}. In Section~\ref{s:nilpotent}, we consider groups with a normal 2-subgroup and prove Lemma~\ref{lem:index2} and Corollary~\ref{cor:nilp}.

\subsection{Preliminaries}
We restrict to finite groups and simple graphs. For a group $G$ with inverse-closed (symmetric) generating set $S\subseteq G \setminus\{1\}$, we let $\Gamma=\Cay(G,S)$ denote the undirected graph with vertices $G$ and edges $\{g, gs\}_{g\in G,s\in S}$. For vertices $v_1, v_2, \dotsc, v_k \in V(\Gamma)$, we let $(v_1,v_2, \dotsc, v_k)$ denote the walk starting at $v_1$ that successively visits $v_2,\dotsc, v_k$. A path is a walk with pairwise distinct vertices; a cycle is a walk $(v_1,v_2, \dotsc, v_k)$ where $k\ge 4$, $v_1=v_k$, and $v_1, v_2,\dotsc, v_{k-1}$ are pairwise distinct. In this note, all cycles are simple cycles. We compose walks $\alpha,\beta$ left to right so that $\alpha \cdot \beta$ describes the walk that first traverses $\alpha$ and then $\beta$. All group actions including translations are taken from the left.

\section{Odd-ordered groups}
\label{s:odd}

In this section, all groups have odd order. We begin by proving Lemma~\ref{lem:main}.

\begin{proof}[Proof of Lemma~\ref{lem:main}]
	As $\abs{G}$ is odd, we may find an inverse-closed set $S_0 \subset S$ in bijective correspondence with $\overline S$ under the quotient $G\to G/A$. For $G_0 = \langle S_0 \rangle \leq G$, the quotient map $G_0 \to \overline G_0:=G_0/(A \cap G_0)$ induces a covering
	\begin{equation}
		\label{eq:cover}
		p: \Cay(G_0,S_0) \to \Cay(\overline G_0, \overline S)
	\end{equation}
	in the topological sense. As $\overline G_0 \iso G/A$, $ \Gamma = \Cay(\overline G_0, \overline S)$ contains a symmetric figure-eight, given as follows. Let $\alpha$  denote the first cycle in $\Gamma$ of power-of-two length, and $\beta = g_0\alpha$, $g_0\in \overline G_0$ denote its translate, up to change of basepoint. Without loss, we may suppose $\alpha$ and $\beta$ start and meet precisely at the identity, so that their corresponding words in $S_0$ reduce to elements $a,b \in A$, respectively. As $\beta = g_0\alpha$, we may decompose $\alpha=\alpha_1\cdot\alpha_2$ where $\alpha_1$ is the path in $\alpha$ that starts at the identity and ends at $g_0^{-1}$ and $\alpha_2$ is the path in $\alpha$ that starts at $g_0^{-1}$ and ends at the identity. Letting $g_1,g_2\in G_0$ denote the reduction of the words $\alpha_1, \alpha_2$, respectively, we have $a=g_1g_2$ and $b=g_2g_1=g_1^{-1}ag_1$.

	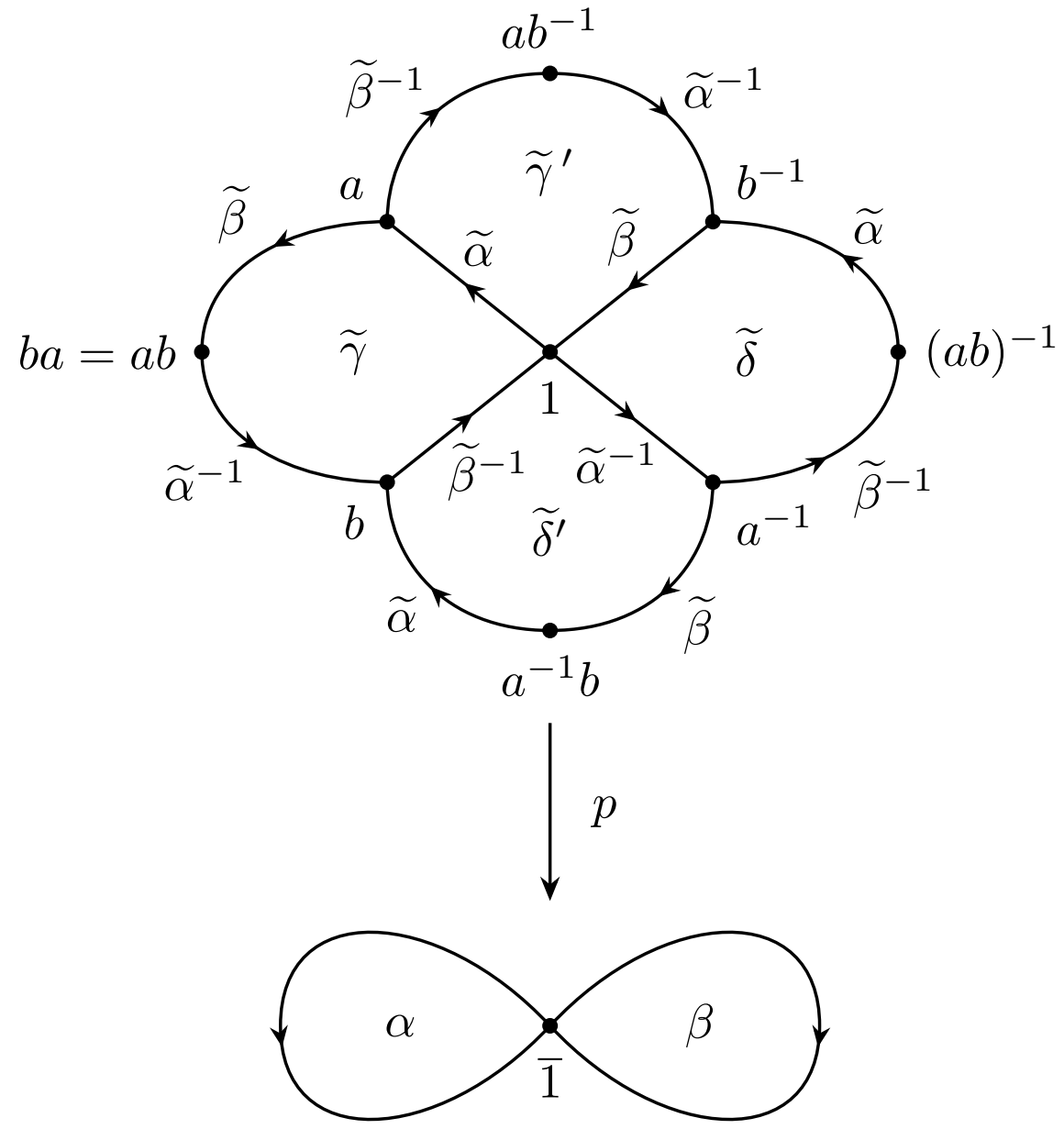
\begin{figure}[ht]
		\centering
		\begin{tikzpicture}[
				x=1.2cm,y=1.2cm,
				line width=.65pt,
				every label/.style={inner sep=3pt},
				path label/.style={inner sep=1pt}
			]

			\begin{scope}[shift={(0,3.35)},xscale=.675,yscale=.54]
				\coordinate (center) at (0,0);
				\coordinate (upper-vertex-left)  at (-1.45,1.45);
				\coordinate (upper-vertex-right) at (1.45,1.45);
				\coordinate (lower-vertex-left)  at (-1.45,-1.45);
				\coordinate (lower-vertex-right) at (1.45,-1.45);
				\coordinate (upper-outer) at (0,3.1);
				\coordinate (lower-outer) at (0,-3.1);

				\draw[path arrow] (center) -- (upper-vertex-left)
				coordinate[pos=.5] (inner-upper-left);
				\draw[path arrow] (lower-vertex-left) -- (center)
				coordinate[pos=.5] (inner-lower-left);
				\draw[path arrow] (upper-vertex-right) -- (center)
				coordinate[pos=.5] (inner-upper-right);
				\draw[path arrow] (center) -- (lower-vertex-right)
				coordinate[pos=.5] (inner-lower-right);

				\foreach \angle/\arrowstyle/\first/\second in {
						0/reverse path arrow/outer-upper-right/outer-lower-right,
						90/path arrow/primed-upper-left/primed-upper-right,
						180/reverse path arrow/outer-lower-left/outer-upper-left,
						270/path arrow/primed-lower-right/primed-lower-left}{
						\begin{scope}[rotate=\angle]
							\draw[\arrowstyle] (1.45,1.45)
							.. controls (2.30,1.45) and (3.1,.95) .. (3.1,0)
							coordinate[pos=.5] (\first);
							\draw[\arrowstyle] (3.1,0)
							.. controls (3.1,-.95) and (2.30,-1.45) .. (1.45,-1.45)
							coordinate[pos=.5] (\second);
						\end{scope}
					}

				\node[path label,above=3pt,xshift=2pt] at (inner-upper-left) {$\widetilde\alpha$};
				\node[path label,below=4pt,xshift=4pt] at (inner-lower-left) {$\widetilde\beta^{-1}$};
				\node[path label,above=3pt,xshift=-2pt] at (inner-upper-right) {$\widetilde\beta$};
				\node[path label,below=4pt,xshift=-3pt] at (inner-lower-right) {$\widetilde\alpha^{-1}$};

				\node[path label,above left=1pt] at (outer-upper-left) {$\widetilde\beta$};
				\node[path label,below left=1pt] at (outer-lower-left) {$\widetilde\alpha^{-1}$};
				\node[path label,above right=1pt] at (outer-upper-right) {$\widetilde\alpha$};
				\node[path label,below right=1pt] at (outer-lower-right) {$\widetilde\beta^{-1}$};

				\node[path label,above left=-1pt] at (primed-upper-left) {$\widetilde\beta^{-1}$};
				\node[path label,above right=1pt] at (primed-upper-right) {$\widetilde\alpha^{-1}$};
				\node[path label,below right=1pt] at (primed-lower-right) {$\widetilde\beta$};
				\node[path label,below left=1pt] at (primed-lower-left) {$\widetilde\alpha$};

				\node[vertex,label={[yshift=-1pt]below:$1$}] at (0,0) {};
				\node[vertex,label={above left:$a$}] at (upper-vertex-left) {};
				\node[vertex,label={below left:$b$}] at (lower-vertex-left) {};
				\node[vertex,label={left:$ba=ab$}] at (-3.1,0) {};
				\node[vertex,label={above right:$b^{-1}$}] at (upper-vertex-right) {};
				\node[vertex,label={below right:$a^{-1}$}] at (lower-vertex-right) {};
				\node[vertex,label={right:$(ab)^{-1}$}] at (3.1,0) {};
				\node[vertex,label={above:$ab^{-1}$}] at (upper-outer) {};
				\node[vertex,label={below:$a^{-1}b$}] at (lower-outer) {};

				\node at (-1.75,0) {$\widetilde\gamma$};
				\node at (1.75,0) {$\widetilde\delta$};
				\node at (0,2) {$\widetilde\gamma\,'$};
				\node at (0,-2) {$\widetilde\delta'$};
			\end{scope}

			\draw[map arrow]
			(0,1.12) -- node[right=5pt] {$p$} (0,.05);

			\begin{scope}[shift={(0,-.7)}]
				\draw[path arrow] (0,0)
				.. controls (-.68,.73) and (-1.62,.77) .. (-1.62,0)
				.. controls (-1.62,-.77) and (-.68,-.73) .. (0,0);
				\draw[path arrow] (0,0)
				.. controls (.68,.73) and (1.62,.77) .. (1.62,0)
				.. controls (1.62,-.77) and (.68,-.73) .. (0,0);
				\node[vertex,label={[yshift=-1pt]below:$\overline 1$}] at (0,0) {};
				\node at (-.9,0) {$\alpha$};
				\node at (.9,0) {$\beta$};
			\end{scope}
		\end{tikzpicture}
		\caption{The cycles $\widetilde\gamma$, $\widetilde\gamma\,'$, $\widetilde\delta$, and $\widetilde\delta '$ in $\widetilde\Gamma$ expressed in terms of lifts of $\alpha$ and $\beta$ in $\Gamma$.}
		\label{fig:gamma_delta}
	\end{figure}

	If $a$ or $b$ were trivial, then both would be, so that the lifts of $\alpha$ and $\beta$ starting from the identity would prescribe a symmetric figure-eight in $\widetilde \Gamma =  \Cay(G_0,S_0)$. If $a=b$, then the lifts of $\alpha \cdot \beta^{-1}$ and $\beta^{-1} \cdot \alpha$ would suffice. The same argument applies for $a=b^{-1}$. We proceed under the hypothesis that $V_0=\{a^{\pm 1},b^{\pm 1}\}$ constitute four distinct points. By the unique lifting property of covering spaces, we need only check for intersections at the endpoints of lifts of $\alpha^{\pm 1}$ and $\beta^{\pm 1}$.

	For $\gamma = \alpha \cdot \beta \cdot \alpha^{-1} \cdot \beta^{-1}$ and $\delta = \alpha^{-1} \cdot \beta^{-1} \cdot \alpha \cdot \beta$, consider their respective lifts $\widetilde\gamma$ and $\widetilde\delta$ starting at the identity in $\widetilde\Gamma$. As $A$ is abelian, $\widetilde \gamma$ and $\widetilde\delta$ are cycles, where the group element $(ab)^{-1}$ translates $\widetilde\gamma$ to $\widetilde \delta$. Let $\gamma' = \alpha \cdot  \beta^{-1} \cdot \alpha^{-1} \cdot  \beta$ and $\delta' =\alpha^{-1} \cdot \beta \cdot \alpha \cdot \beta^{-1}$ correspond to reversing $\beta$ in $\gamma$ and $\delta$. Then $\widetilde \gamma\, ' $ and $\widetilde \delta' $, defined analogously, are also cycles and translates of each other. See Figure~\ref{fig:gamma_delta} for a schematic.

	We claim that $\widetilde\gamma, \widetilde\delta$ or $\widetilde\gamma', \widetilde\delta'$ constitute a symmetric figure-eight in $\widetilde\Gamma$.  As $ab$ and $ab^{-1}$ are nontrivial and distinct from their inverses, it remains to show that $ab \notin V_0$ or $ab^{-1} \notin V_0$. Suppose $ab=a^{-1}$. Then $b=a^{-2}$ so that the distinctness of $V_0=\{a^{\pm 1},a^{\pm 2}\}$ implies $\abs{a} \ge 5$.  If $\abs{a}=5$, the identity $b=g_1^{-1}ag_1$ becomes $a^3=g_1^{-1}ag_1$, so that conjugation by $g_1$ on $a$ would induce an automorphism of $\langle a\rangle$ of order four, contradicting the odd order of $g_1$. Consequently, $\abs{a} \ge 7$ so that $ab^{-1} = a^3 \notin V_0$. An analogous argument applies for the case $ab=b^{-1}$.
\end{proof}

We now consider Theorem~\ref{thm:main}. Lemma~\ref{lem:main} allows us to argue that a minimal counterexample would take the form of a Frobenius group. We then handle the resulting cases.

\begin{proof}[Proof of Theorem~\ref{thm:main}.]
	For a cyclic group $G$ with inverse-closed generating set $S \subset G \setminus \{1\}, \abs{S}\ge 3$, let $s,t\in S$ such that $s\ne t^{\pm 1}$. Then $s$ and $t$ commute, so that the closed walk of length 4 corresponding to the word $[s,t]=sts^{-1}t^{-1}$ is a cycle.

	We claim that for all noncyclic groups $G$ of odd order, $\Cay(G,S)$ contains a symmetric figure-eight for any inverse-closed generating set $S$, $\abs{S} \ge 3$.
	Suppose by way of contradiction that $\Gamma=\Cay(G,S)$ is a counterexample with $\abs{G}$ minimal. If $G$ were abelian, then $(1,s,st,t,1)$ and its translate by $st$ would give a symmetric figure-eight for any $s,t \in S, t \notin \langle s \rangle$. If $G_0 = \langle s,t\rangle \lneq G$ for noncommuting $s,t\in S$, then $\Cay(G_0, \{s^{\pm 1},t^{\pm 1}\})$ would also be a counterexample, contradicting the minimality of $\abs{G}$. We proceed under the assumption that $G$ is nonabelian and $S=\{s^{\pm 1},t^{\pm 1}\}$.

	By Feit--Thompson, $G$ is solvable. If $G/N$ were noncyclic for some nontrivial abelian $N\triangleleft G$, the minimality of $G$ and Lemma~\ref{lem:main} would allow us to conclude. Let $V\triangleleft G$ be minimal normal and thus elementary abelian. We must have $G/V$ cyclic so that $[G,G] = V$. Additionally, $Z( G)=1$, as otherwise $ G/Z( G)$ would be cyclic, making $G$ abelian. Find $x\in G$ such that $x V$ generates $ G/V$. Let $m$ denote the order of $xV$ in $G/V$. As $x^m \in V$ lies in $Z(G)$, it follows that $x^m=1$ so that $ G = V \rtimes \langle x \rangle$. As $V$ is minimal normal, the action of $\langle x \rangle$ on $V$ is irreducible, and $Z( G)=1$ implies $x$ acts fixed-point-freely. For each $1<i < m$, $Z_V(x^i)$ is $x$-invariant, and so either trivial or all of $V$. However, $Z_V(x^i)=V$ would imply $x^i \in Z(G)$, forcing $Z_V(x^i)$ to be trivial. Consequently, $G$ is Frobenius.

	We identify $V$ with the additive group of a finite field $F$ of characteristic $p$ so that $x$ acts by multiplication by some $r\in F^\times$, i.e. $xvx^{-1}=rv$. Let $\mu_r$ denote the minimal polynomial of $r$ over $\FF_p$, equivalently that of conjugation by $x$ on $V$. Irreducibility implies $F=\FF_p(r)\iso \FF_p[X]/(\mu_r)$ and $p \nmid m$ since $\abs{r} = m$ divides $\abs{F}-1$. Using the terminology of affine groups~\cite[\S 2.8]{Dix96}, we let $\tau_c(z) = z + c$ denote translation by $c\in F$ so that $x\tau_c x^{-1} = \tau_{rc}$. In particular, every element of $G$ may be expressed uniquely as $\tau_a x^i$, where $a\in F$ and $0\leq i < m$. For $a,b \in F$ and $0 \leq i,j < m$, we have  $(\tau_a x^i)(\tau_b x^j) = \tau_{a + r^{i} b} x^{i+j}$ where $i+j$ is taken modulo $m$.

	Let $\overline s, \overline t$ be the images of $s, t$, respectively, in $\overline G =G/V$.  For any automorphism $\sigma \in \operatorname{Aut}(G)$, we have $\Cay(G, S) \iso \Cay(G , \{\sigma(s): s \in S\})$. We make use of this property  to select generators in the cases that follow. In cases (2)--(4), we may suppose
$\overline s$ generates $G/V$, after exchanging generators if necessary. If $\overline s \ne 1$, then $s$ has a unique fixed point and conjugation by a translation sends that point to $0$. Since $\overline s$ generates $G/V$, the conjugated element is a generator of the scalar subgroup. Without loss, we may set $x=s$ and denote its scalar by $r$. Then, $t=\tau_a x^i$ where $a\ne 0$ since $s$ and $t$ do not commute and $i=2,1,0$ respectively, after inverting $t$ if necessary. Conjugation by the dilation $d_a(z)=a^{-1}z$ normalizes $G$, fixes each $x^i$, and sends $\tau_a x^i$ to $\tau_1x^i$.

	\begin{enumerate}
		\item \emph{Suppose $\overline s, \overline t$ and their inverses are distinct and nontrivial, and $m\ne 5$.} After replacing $\overline t$ with $\overline t{}^{-1}$ if necessary, we have $\overline s \ne \overline t{}^{-2}$ and $\overline t \ne \overline s{}^{-2}$, for if neither $\overline t$ nor $\overline t{}^{-1}$ satisfied both equations, then $\overline G=\langle\overline s,\overline t\rangle
\iso C_5$, contradicting $m\ne5$. The cycle $\alpha=(1, \overline s, \overline s\overline t, \overline t, 1)$ and its translate $\overline s\overline t\alpha$ then form a symmetric figure-eight in $\Cay(\overline G, \{\overline s^{\pm 1}, \overline t^{\pm 1}\})$. Lemma~\ref{lem:main} allows us to conclude.

		\item \emph{Suppose $\overline s, \overline t$ and their inverses are distinct and nontrivial, and $m = 5$.} Without loss, we may suppose $s=x$ and $t=\tau_1x^2$, so that $[s,t], [t,s^{-1}] \in V$. Then $[[s,t], [t,s^{-1}]] = sts^{-2}t^{-1}stst^{-1}s^{-2}tst^{-1} = sts^3t^{-1}stst^{-1}s^3tst^{-1}$ is a word of length 16, where the last equality follows from $s^{-2}=s^3$. The walk $\beta$ corresponding to this word traces a simple cycle, with vertices $(1,x,\tau_r x^3,\tau_r x^4,\tau_r,\tau_r x,\tau_{r-r^4}x^4,\tau_{r-r^4},\tau_{1+r-r^4}x^2,\tau_{1+r-r^4}x^3,\tau_{1-r^4}x,$ \\ $\tau_{1-r^4}x^2,\tau_{1-r^4}x^3,\tau_{1-r^4}x^4,\tau_1x,\tau_1x^2,1)$. Likewise, the walk $\beta'$ corresponding to replacing $s$ and $t$ with $t$ and $s^{-1}$, respectively, in the definition of $\beta$ is a simple cycle of length 16. Then, $\beta$ and $\tau_1\beta$ give a symmetric figure-eight with common vertex $\tau_1x$ unless $r^4 \in \{r-1,r+1,2\}$ or $r = 2$; i.e. when $(p,r) \in \{(11,3),(11,4),(31,2),(31,16)\}$.  In these exceptional cases, $\beta'$ and $\tau_{-(1+r^2)}\beta'$ give a symmetric figure-eight with common vertex $\tau_{-r^2}x^2$.

		\item \emph{Suppose $\overline s = \overline t \ne 1$ in $\overline G$, after inverting a generator if necessary.} Without loss, $s=x$ and $t=\tau_1x$. Then $ts^{-1}, s^{-1}t \in V$ and we can consider the cycle $\gamma$ corresponding to the commutator $[ts^{-1},s^{-1}t] = ts^{-2}tst^{-2}s$ with vertices $(1,\tau_1x,\tau_1,\tau_1x^{-1},\tau_{1+r^{-1}},\tau_{1+r^{-1}}x,\tau_{r^{-1}},x^{-1},1)$. It will intersect at precisely $\tau_{1+r^{-1}}$ with $\tau_{1+r^{-1}}\gamma$ if $r^{-1}\notin \{-2,-1/2\}$ and at precisely $\tau_1$ with $\tau_{1-r^{-1}}\gamma$ if $r^{-1}\notin \{2,1/2\}$. For both translates to fail to form a symmetric figure-eight, it would need to be the case that $r^{-1} \in \{-2,-1/2\}\cap\{2,1/2\}$. As $p$ is odd, it would follow that $p=5$ and $r \in \FF_5^\times \iso C_4$, contradicting that $m$ is odd.

		\item \emph{Suppose either $\overline s =1$ or $\overline t =1$.} Without loss, the generators can be taken to be $s=x$ and $t=\tau_1$. As $t$ and $sts^{-1}=\tau_r$ commute, the cycle $\delta$ corresponding to the commutator $[t,sts^{-1}]$ is a simple eight-cycle with vertices $(1,\tau_1,\tau_1x,\tau_{1+r}x,\tau_{1+r},\tau_r,\tau_r x,x,1)$. It intersects  $x\tau_1\delta$ at precisely $x\tau_1$ if $r-r^{-1} \notin \{-1,-2,-3/2\}$ or with $x\tau_{-1}\delta$ at precisely $x$ if $r-r^{-1}\notin \{1,2,3/2\}$. From this, we need only consider $p\leq 7$. As $r-r^{-1} \in \FF_p$, $f(X) =X^2-(r-r^{-1})X- 1\in \FF_p[X]$ has root $r$ so that $\mu_r \mid f$. If $\mu_r$ had degree two, then $f=\mu_r$ so that $f\mid (X^m -1)$, forcing both roots of $f$, namely $r$ and $-r^{-1}$, to be $m$-th roots of unity, a contradiction. Consequently, $\mu_r$ has degree one; i.e. $r\in \FF_p$ so that $F=\FF_p$. Characteristics three and five admit no nontrivial scalars of odd order $m>1$, leaving $G\iso C_7 \rtimes C_3$ as the sole exception. Reversing $x$ if need be, let $sts^{-1}= t^2$. The cycle $\delta'=(1,\tau_1,\tau_2,\tau_3,\tau_4,\tau_5,\tau_5x,x,1)$ corresponding to the word $t^5sts^{-1}$ and its translate $x\delta'$ provide a symmetric figure-eight in this case.
	\end{enumerate}

	We find that $\Gamma$ contains a symmetric figure-eight in every case, allowing us to conclude.
\end{proof}

\section{Finite nilpotent groups}
\label{s:nilpotent}

We begin with the proof of Lemma~\ref{lem:index2}.

\begin{proof}[Proof of Lemma~\ref{lem:index2}]
	Without loss, $P\lneq G$ and $S\cap P$ contains no element of order at least four. If $s \in S\cap P$ is an involution, let $t\in S\setminus P$ and consider $x=[s,t^{-1}]=st^{-1}st \in P$. As $\abs{x}$ is a power of two, the walk $\alpha$ corresponding to the word $(st^{-1}st)^{\abs{x}}$ is a cycle of length $4\abs{x}$. Its vertices lie in the four  sets $\langle x \rangle,  \langle x \rangle s,  \langle x \rangle st^{-1}, \langle x \rangle t^{-1}$. The first two lie in $P$ and the last two in $t^{-1}P$. If $s\in\langle x \rangle$, then $sx=t^{-1}st \in \langle x \rangle$ as well, forcing $s=sx$ to be the unique nonidentity involution in $\langle x \rangle$, a contradiction. Thus, $\langle x \rangle \ne  \langle x \rangle s$ so that $\langle x \rangle st^{-1} \ne \langle x \rangle t^{-1}$, making $\alpha$ a simple cycle. We proceed under the hypothesis $S \cap P = \emptyset$.

	Suppose $\overline s = \overline t$ for distinct $s,t \in S$. Then $y=s^{-1}t \in P$ so $\abs{y}$ is a power of two. Let $\beta$ be the walk associated to the word $(s^{-1}t)^{\abs y}$. Its vertices alternate between the disjoint sets $\langle y \rangle$ and $\langle y \rangle t^{-1}$, so that it is a simple cycle of length $2\abs{y}$.

	Otherwise, $S \to \overline S$ is bijective so that the quotient $G \to G/P$ induces a covering
	\begin{equation}
		p: \Cay(G,S) \to \Cay(G/P, \overline S)
	\end{equation}
	in the topological sense. Let $\gamma \in \Cay(G/P, \overline S)$ denote a simple cycle of length $2^n, n\ge 2$. Without loss, we may suppose it starts at the identity $\overline 1 \in G/P$. Each vertex in the connected component of $p^{-1}(\gamma)$ containing $1 \in G$ will have degree 2. Let $\widetilde\gamma$ denote this cycle, starting at $1$. If $w=s_1s_2\cdots s_{2^n}$ is the word in $S$ such that $\gamma$ traces $\overline s_1 \overline s_2\cdots \overline s_{2^n}$, then $\widetilde\gamma$ will have length $2^n \abs{w}$, where $w \in P$.
\end{proof}

Corollary~\ref{cor:nilp} now follows.

\begin{proof}[Proof of Corollary~\ref{cor:nilp}]
	By Theorem~\ref{thm:main}, we may suppose $\abs{G}$ is even. Let $P \triangleleft G$ be its Sylow 2-subgroup. If $G=P$ and $S$ consists entirely of involutions, then for distinct $s,t \in S$, $(st)^{\abs{st}}$ traces a simple cycle of length $2\abs{st}$; if $G=P$ and $s\in S$ is not an involution, the word given by $s^{\abs{s}}$ is a cycle of power-of-two length. Thus, $P \lneq G$ is proper. If $S\cap P \ne\emptyset$ or two distinct elements of $S$ map to the same image in $G/P$, the first two paragraphs of Lemma~\ref{lem:index2} allow us to conclude. Otherwise, $S\to \overline S = \{sP\}_{s\in S}\setminus \{P\}$ is bijective and Theorem~\ref{thm:main} implies $\Cay(G/P, \overline S)$ contains a power-of-two-length cycle, allowing us again to appeal to Lemma~\ref{lem:index2}.
\end{proof}

\subsection*{Acknowledgements}
The author is grateful to the University of Chicago Library for providing scanned copies of cited proceedings papers from Congressus Numerantium.


\end{document}